\documentclass{article}
\usepackage{amssymb}

\usepackage{graphicx}
\usepackage{amsmath}

\newtheorem{theorem}{Theorem}

\newtheorem{lemma}[theorem]{Lemma}

\newtheorem{proposition}[theorem]{Proposition}

\begin{document}

\title{Varadhan estimates for a degenerated convolution semi-group: upper bound}
\author{R\'emi L\'eandre   \\
Institut de Math\'ematiques, Universit\'e de Bourgogne,
\\ 21000, Dijon, FRANCE.
} \maketitle

\begin{abstract}
We give a proof of Varadhan estimates for a denerated jump process with independent increments with more and more jumps which become smaller and smaller.
 The proof uses the Malliavin Calculus of Bismut type for jump process in semi-group theory of [14] and [15]
and 
Wentzel-Freidlin estimates for jump processes in semi-group theory of [16] and [17].
\end{abstract} 
\section{Introduction}
The object of  the large deviation theory [22] is to estimate the logarithm of the probability of rare events. Bismut [1] pointed out the relationship
 between large deviation estimates and the Malliavin Calculus in order to get short time asymptotics of heat kernels associated to diffusion semi-groups.
This relationship was fully performed by L\'eandre in [8], [9].
 The reader interested in short time asymptotics of heat-kernels by using the Malliavin Calculus as a tool
 can look at the 
review of L\'eandre [10], Kusuoka [5] and Watanabe [21].

L\'eandre  has translated plenty of tools of stochastic analysis in semi-group theory  by using the fact that in the proof of these tools,
 there are suitable stochastic
 differential equations and therefore suitable parabolic equations which appear. We refer to the review of L\'eandre [11], [13] on that.

In [6], [7], we established the existence of a density for a degenerated convolution semi-group by using
 the Malliavin Calculus of Bismut type as a tool.
We have translated the proof of our result in semi-group theory in [15].

In [16], [17], we have translated the proof of Wentzel-Freidlin estimates for jump processes in semi-group theory [22]. The goal of this paper is to do the marriage
between Wentzel-Freidlin estimates  in semi-group theory and the Malliavin Calculus of Bismut type for jump process in semi-group theory [14] in order to get
 Varadhan estimates for a degenerated convolution semi-group. Let us recall that this marriage was done first in the stochastic case
 by Ishikawa [3] and 
Ishikawa-L\'eandre [4].

Let us recall that there are classically 3 types of estimates for heat-kernels:

-)Rough but global estimates of Aronson type, which take care of the Gaussian decay outside the diagonal as well as their polynomial blowing-up on 
the diagonal. These estimates reflect the global properties of the generator [2].

-)Precise and local estimates of Pleijel-Minakishudaram type [19]. These precise estimates are required too in
 the theory of semi-classical expansion of the solution of the Schroedinger equation [18]. It is still these precise asymptotics
 which play a role in index theory and spectral asymptotics [?].

-)Varadhan type estimates [20] which say when Wentzel-Freidlin estimates pass to heat-kernels.

Our estimate works in the case of a generator which is not a pseudo-differential operator: our result presents some difficulties to be reached by the classical technics of microlocal analysis.

\section{Statement of the main theorems}
Let $g(z)$ be a smooth function with compact support in $[-1,1]$  from
$\mathbb{R}^*$ into $\mathbb{R}^+$ such that
\begin{equation}\int_{\mathbb{R}}(z^2\wedge 1)g_j(z)dz < \infty\end{equation}
Moreover we suppose that on a neighborhood of 0
\begin{equation}g(z) = {C \over \vert z\vert^{1+\alpha}}\end{equation}
for some $\alpha \in ]0,1[$. 

We consider a smooth curve  $\gamma(z)$ into $\mathbb{R}^d$ with bounded derivatives of each order such that
\begin{equation}\gamma(0) = 0\end{equation}
We do the following hypothesis:

 {\bf{Hypothesis H:}}There exists a $k$ such that
$\cup_{ l\leq k}{d^l\over dz^l}\gamma(0)$ spans 
$ \mathbb{R}^d$.

We consider the Markov generator acting on $C_b^\infty(\mathbb{R}^d)$
\begin{multline}Lf(x) =  \int_{{R}}(f(x+\gamma(z))-f(x))
g(z)dz = \\ \int (f(x+z')-f(x))\mu(dz')\end{multline}
$\mu$ is the Levy measure associated to the jump process associated to $L$.
 We introduce the Hamiltonian defined for $ \xi \in  \mathbb{R}^d$
\begin{equation}H(x, \xi) = \int_{\mathbb{R}^d}(\exp[<z',\xi>]-1)\mu(dz')\end{equation}
We consider the Legendre transform of $H$
\begin{equation}L( \alpha) = \sup_{\xi}\{<\alpha,\xi>-H(\xi)\}\end{equation}
We will show later that the convex function $L$ is finite (therefore $L$ is continuous).
 
We consider a piecewise $C^1$ curve $\phi(t)$ and we consider the action:
\begin{equation} S(\phi) = \int_0^1L( d/dt\phi(t))dt\end{equation}
We put
\begin{equation}l(x,y) = \inf_{\phi(0)=x, \phi(1) = y}S(\phi) \end{equation}
Under the previous assumption, $(x,y) \rightarrow l(x,y)$ is continuous.

$L$ spans a convolution Markovian semi-group $P_t$.
Let us recall [15] that $P_t$ has a smooth heat-kernel:
\begin{equation}P_tf(x) =\int_{\mathbb{R}^d}p_t(x,y)f(y)dy\end{equation}

We define the generator $L^h$ defined on smooth functions with bounded derivatives at each order
\begin{equation} L^hf(x) = \\ \int_{\mathbb{R}^d}(f(x+hz')-f(x))\mu(dz')\end{equation}

Under these previous assumptions, we will get in the sequel Markovian semi-group. In particular, $1/hL^h$ generates a semi-group $P_t^h$.
In such a case, $P_t^h$ is a semi-group in probability measures ([22]). We get large deviation results of Wentzel-Freidlin type for this semi-group:

\begin{theorem} Let $O$  be an open ball relatively compact of $\mathbb{R}^d$. When $h \rightarrow 0$,
\begin{equation}\overline{\lim} h\log P_1^h[1_O](x) \leq -\inf_{y \in O}l(x,y)\end{equation}
\end{theorem}
Under the previous assumptions, the previous estimate passes to heat-kernel. This means that the following Varadhan type estimate holds:
\begin{theorem} When $h \rightarrow 0$, we have uniformly on any compact of $\mathbb{R}^d$
\begin{equation}\overline{\lim} h\log p_1^h(x,y)\leq  -l(x,y)\end{equation}
\end{theorem}

\section{Proof of the theorems}
In order to simplify the proof, we will suppose that ${d^l \over dz^l}\gamma(0)$ are independent, $l \leq d$
\begin{lemma}:For any $K$, there exists $C_K$ such that for $\vert \alpha \vert \geq C_K$, $L(\alpha)\geq K \vert \alpha\vert$.
 Moreover $L$ is continuous, convex and tends to 
$\infty$ when $\vert \alpha \vert \rightarrow \infty$\end{lemma}
{\bf{Proof:}} We work in the basis given by $d^l/dz^l \gamma(0)$. We get if $\xi_l \geq 0$
\begin{multline}H(\xi) \geq \sum_lC_l\int_0^1(\exp[C_l<\xi_l,z^l>]-1)z^{-1-\alpha}dz+\\  \sum_lC_2\int_0^1(\exp[-C_2<\xi_l,z^l>]-1)z^{-1-\alpha}dz
\end{multline}
In each integral, we perform the change of variable $z' = z^l$. We get a lower bound of the type
\begin{multline}
H(\xi) \geq \sum_lC_l\int_0^1(\exp[C_l<\xi_l,z>]-1)z^{-1-\alpha/l}dz + \\ \sum_lC_2\int_0^1(\exp[-C_2<\xi_l,z>]-1)z^{-1-\alpha/l}dz\end{multline}
By putting $z' = \xi_l z$, we get the bound
\begin{equation}H(\xi) \geq \sum_lC_l\exp[C_l\xi_l]\end{equation}
Analogous considerations work in the others sectors of $\mathbb{R}^d$.
Therefore
\begin{equation}L(\alpha) \leq \sup_\xi ( \sum \alpha_l\xi_l-\sum C_l\exp[C_l\xi_l])\leq C\vert\alpha\vert \log (\vert \alpha \vert)+C\end{equation}
We have clearly the opposite inequality of (14) and therefore
\begin{equation}L(\alpha) \geq C\vert\alpha\vert \log (\vert \alpha \vert)\end{equation}

 By standard argument, $L$ is convex.
Since $L$ is finite, it is therefore continuous.$\diamondsuit$.

{\bf{Proof of Theorem 1}}
It it is an adaptation of the proof of the analogous theorem in [17]. Since $L$ is continuous, we cand find $\alpha_i$ and $\xi_i$ such that if
$\vert \alpha \vert < R$ and if we put 
\begin{equation}L'(\alpha) = \sup_i\{L(\alpha_i) + <\xi_i, \alpha-\alpha_i>\}\end{equation}
then
\begin{equation}L(\alpha)-L'(\alpha) \leq \chi\end{equation}
for a small $\chi$.

But the Legendre transform of $L$ is $H$. Therefore we can find the $\xi_i$ such that
$-L(\alpha_i)+<\xi_i, \alpha_i>$ are close to $H(\xi_i)$. In order to be more precise,
 we can repeat the considerations of the proof of the same theorem in [17] by using the fact that
 a continuous convex function is differentiable on a dense set.
$\diamondsuit$

 We consider a smooth function $\nu(z)$ with compact support with values in 
$\mathbb{R}^*$ equal to $z^4$ in a neighborhood of 0.
We consider the space $\mathbb{R}^d\times \mathbb{M}_d$ where $\mathbb{M}_d$ is the space of 
symmetric  matrices on $\mathbb{R}^d$. $(x,V) \in \mathbb{R}^d \times\mathbb{ M}_d$. We consider the Malliavin generator
 acting on test functions $\hat{f}$ on this space
\begin{equation}\hat{L}^h\hat{f}(x,V) =  \int_{\mathbb{R}}(\hat{f}(x+\gamma(hz),
V+\nu(hz)<.,\gamma'(hz)>^2)-\hat{f}(x,V))g(z)dz
\end{equation}
$V$ is called the Malliavin matrix.
$1/h\hat{L}^h$ generates a probability semi-group.
\begin{proposition}For all $p$, there exists $n(p)$ such that
\begin{equation}\hat{P}_t^h[V^{-p}](x,0) \leq C h^{-n(p)}\end{equation}
for $h \leq 1$.
\end{proposition}
{\bf{Proof:}}
Let us recall (Lemma 6 of [?]) that for all $p>0$, all 
\begin{equation}\sup_{h\leq 1}\hat{P}^h_t[V^p](x,0) < \infty \end{equation} 
such that is enough to show for all $p > 0$
\begin{equation}\sup_{\vert \xi\vert = 1}\hat{P}^h_t[<V,\xi>^{-p}] < Ch^{-n'(p)} \end{equation}
We consider the semi-group $P_t^{\xi,h}$ associated to the generator $L^{\xi,h}$ on 
$\mathbb{R}^d \times \mathbb{R}$, $(x,u) \in \mathbb{R}^d \times \mathbb{R}$:
\begin{equation}L^{\xi,h}\hat{f}(x,u) = 1/h \int_{\mathbb{R}}(\hat{f}(x+\gamma(hz), u+\nu(hz)<\xi,\gamma(hz)>^2)-\hat{f}(x,u))
g(z)dz
\end{equation}
We have as it was remarked in [15]
\begin{multline}\hat{P}^h_t[<V,\xi>^{-p}](x,0) = P_t^{\xi,h}[u^{-p}](x,0) \\=
\Gamma(p)^{-1}\int_0^\infty\beta^{p-1}P_t^{\xi,h}[\exp[-\beta u]](x,0)d\beta\end{multline}
We consider the generator $L^{\beta,h}$ acting on functions on $\mathbb{R}^d\times \mathbb{R}\times \mathbb{R}$ ($(x,u,v) \in \mathbb{R}^d \times 
\mathbb{R}\times \mathbb{R})$
\begin{multline}L^{\beta,h} \hat{f}(x,u,v) =\\ 1/h \int_{\mathbb{R}}(\hat{f}(x+\gamma(hz), u+<\xi,\gamma'(hz)>^2\nu(hz),v)- \hat{f}(x,u,v))g(z)dz\\
+ 1/h\int_{\mathbb{R}}<\exp[-2\beta\nu(hz)<\xi,\gamma'(hz)>^2]-1,grad_v\hat{f}(x,u,v)>g(z)dz\end{multline}
Let us consider the function $\exp[-2\beta u - v]$. If we apply $L^\beta$ to it, we find zero: Therefore
\begin{equation}P_1^\beta[\exp[-2\beta u - v]](x,0,0) = 1\end{equation}
Since we have Markov semigroups, we get the following lemma [15] originally proved by Bismut by using martingales theory:
\begin{lemma}(Bismut[?]
)For all positive $p$
\begin{equation}\hat{P}^h_t[<V,\xi>^{-p}](x,0) \leq\Gamma(p)^{-1}\int_0^\infty\beta^{p-1}(P_t^{\beta,h}[\exp[v]](x,0,0))^{1/2}d \beta
\end{equation}
\end{lemma}
But 
\begin{multline}P_t^{\beta,h}[\exp[v]](x,0,0) = \\ \exp[1/h t\int_{\mathbb{R}}(\exp[-2\beta\nu(hz)<\gamma'(hz),\xi>^2]-1)g(z)dz]\end{multline}
and we deduce our results by the  uniform Tauberian theorem of [5], [6] and [15].
Namely
\begin{multline}1/h t\int_{\mathbb{R}}(\exp[-2\beta\nu(hz)<\gamma'(hz),\xi>^2]-1)g(z)dz \leq \\
1/ht\int_1^1(\exp[-C\beta h^4z^4\sum(\xi_l,( hz)^{l-1})^2]-1)-1)g(z)dz\end{multline}
the result holds by the uniform Tauberian theorems of the Annexe.$\diamondsuit$

We put $\mathbb{E}^l = \mathbb{R}^d \times \mathbb{M}_d\times \mathbb{R}^{d_2}\times..\times \mathbb{R}^{d_l}$, $x^l = (x_1, v, x_2,..,x_l)$
 is the generic element of this big space .
 We consider a bounded map $G(z)$ from $\mathbb{R}$ into $\mathbb{E}^l$. Its values in $\mathbb{R}^d$ are $\gamma(z)$. Its values in $\mathbb{M}_d$ are
$\nu(z)<\gamma'(z),.>^2$. 
 The other components have compact support and are bounded by $C\vert z\vert^2$.

We consider the generator

\begin{equation}L^{l,h}f^l(x^l) =
 1/h\int_{\mathbb{R}}(f^l(x^l+ G(hz))\\-f^l(x^l))g(z)dz\end{equation}
It generates a Markov semi-group $P_t^{l,h}$. 

The main remark whose proof follows the same lines as the proof of Theorem 6 of [14] is the following proposition:
\begin{lemma}For a function $f^l$ with polynomial growth
\begin{equation}\sup_{h\leq 1}P_t^{l,h}[\vert f^l\vert](x^l) < \infty
\end{equation}\end{lemma}
{\bf{Proof of Theorem 2:}}
We consider the Markov generator
\begin{multline}
L^{1,h}f(x_1,V) = \\1/h\int_{\mathbb{R}}(f(x_1 + \gamma(z), V +\sum \gamma'(z)\nu(z)<\gamma'(z),.>)-f(x_1,V))g(z)dz\end{multline}
We enlarge the semi-group given by this Markov generator by the system given before.
We get the integration by parts of Bismut type to this test function in order to get for all multi-index 
$(\alpha) = (\alpha^1,.., \alpha_d) \in \mathbb{N}^d$. $D^{(\alpha)} = {\partial^{\sum \alpha_i}
\over\partial x_1^{(\alpha_1)}..\partial x_d^{(\alpha_d)}}$
\begin{equation}P_t^h[D^{(\alpha)}f](x) = 
P_t^{l,h}[f\theta](x,0)\end{equation}
where $P_t^{l,h}$ is a semi-group of the type considered before and $\theta$
 a polynomial in the extra-variables  and in the inverse of the Malliavin matrix. We consider a smooth
 function $g$ positive equal to 1 at $y$ with a small support. We apply the previous inequality to the the test function $fg$.
We deduce by Theorem 1, Proposition 4 , Lemma 6 and by induction on the length of $(\alpha)$
that
\begin{equation}
\vert P_1^l[gD^{\alpha}f](x)\vert \leq Ch^{-n(\alpha)}\exp[{-l(x,y)+\chi \over h}]\end{equation}
for a small $\chi$. This shows the result.$\diamondsuit$  
\par\noindent

\section{Annexe: Tauberian theorems}
Let us recall the uniform Tauberian theorem of [6], [7], [15].

Let $F$ be a set of smooth maps $f$ from $\mathbb{R}$ into $\mathbb{R}$ equal to 0 in 0 such for all $k$
 \begin{equation}\sup_{f \in F, z \in \mathbb{R}}\vert f^{(k)}(z)\vert < \infty \end{equation}
and such that there exists $K$
such that
\begin{equation}\inf_{f\in F}(\sup_{k\leq K}\vert f^{(k)}(0)\vert)> C > 0\end{equation}

Let us define
\begin{equation}\tau_f(\beta) = \int_{\mathbb{R}}(\exp[-\beta\vert f(z)\vert]-1)g(z)dz\end{equation}
where $\int_{\mathbb{R}}(z^2\wedge 1)g(z)dz < \infty$ and where $g(z) \geq 0$
Let us suppose that for some $\alpha \in ]0,2[$
\begin{equation}\underline{\lim}_{z \rightarrow 0}\vert z\vert^{1+\alpha}g(z) > 0\end{equation}
Then
\begin{theorem}If (36), (37), and (39) are checked, there exists $\alpha_1 > 0$ which depends only from the data in (36), (37), (38) and (39) such that
\begin{equation}\overline{\lim}_{\beta \rightarrow \infty}(\sup_{f \in F}{\tau_f(\beta)\over \beta^{\alpha_1}}) < 0\end{equation}\end{theorem}

\end{document}